\documentclass{IEEEtran4PSCC}
\ifCLASSINFOpdf
   \usepackage[pdftex]{graphicx}
\else
   \usepackage[dvips]{graphicx}
\fi
\usepackage{amsmath}
\usepackage{booktabs}
\usepackage{mathtools}
\usepackage{graphicx}
\usepackage{xcolor}
\usepackage[caption=false]{subfig}

\makeatletter
\let\old@ps@headings\ps@headings
\let\old@ps@IEEEtitlepagestyle\ps@IEEEtitlepagestyle
\def\psccfooter#1{%
    \def\ps@headings{%
        \old@ps@headings%
        \def\@oddfoot{\strut\hfill#1\hfill\strut}%
        \def\@evenfoot{\strut\hfill#1\hfill\strut}%
    }%
    \def\ps@IEEEtitlepagestyle{%
        \old@ps@IEEEtitlepagestyle%
        \def\@oddfoot{\strut\hfill#1\hfill\strut}%
        \def\@evenfoot{\strut\hfill#1\hfill\strut}%
    }%
    \ps@headings%
}
\makeatother

\psccfooter{%
        \parbox{\textwidth}{\hrulefill \\ \small{22nd Power Systems Computation Conference} \hfill \begin{minipage}{0.2\textwidth}\centering \vspace*{4pt} \includegraphics[scale=0.06]{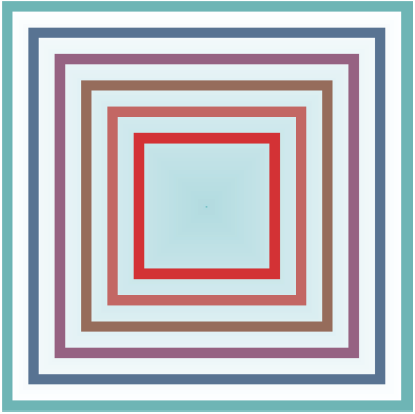}\\\small{PSCC 2022} \end{minipage} \hfill \small{Porto, Portugal --- June 27 -- July 1, 2022}}%
}

\begin{document}
%
\title{Optimal Power Flow in Four-Wire Distribution Networks: Formulation and Benchmarking}

\author{
\IEEEauthorblockN{Sander Claeys\\ Geert Deconinck}
\IEEEauthorblockA{KU Leuven / EnergyVille\\
Leuven / Genk, Belgium\\
\{sander.claeys, geert.deconinck\}@kuleuven.be}
\and
\IEEEauthorblockN{Frederik Geth}
\IEEEauthorblockA{CSIRO Energy\\
Newcastle, Australia\\
frederik.geth@gmail.com
}
}


\maketitle

\begin{abstract}
In recent years, several applications have been proposed in the context of distribution networks.
Many of these can be formulated as an optimal power flow problem, a mathematical optimization program which includes a model of the steady-state physics of the electricity network.
If the network loading is balanced and the lines are transposed, the network model can be simplified to a single-phase equivalent model.
However, these assumptions do not apply to low-voltage distribution networks, so the network model should model the effects of phase unbalance correctly.
In many parts of the world, the low-voltage distribution network has four conductors, i.e. three phases and a neutral.
This paper develops OPF formulations for such networks, including transformers, shunts and voltage-dependent loads, in two variable spaces, i.e. current-voltage and power-voltage, and compares them for robustness and scalability.
A case study across 128 low-voltage networks also quantifies the modelling error introduced by Kron reductions and its impact on the solve time. 
This work highlights the advantages of formulations  in current-voltage variables over power-voltage, for four-wire networks.
\end{abstract}


\thanksto{\noindent Submitted to the 22nd Power Systems Computation Conference (PSCC 2022).}

\newcommand{\TODO}[1]{\begin{color}{blue}#1\end{color}}

\newcommand{\cscal}[1]{\textcolor{blue}{#1}}
\newcommand{\rscal}[1]{\textcolor{red}{#1}}
\newcommand{\cvec}[1]{\textcolor{blue}{\mathbf{#1}}}
\newcommand{\rvec}[1]{\textcolor{red}{\mathbf{#1}}}
\newcommand{\cmat}[1]{\textcolor{blue}{\mathbf{#1}}}
\newcommand{\rmat}[1]{\textcolor{red}{\mathbf{#1}}}

\newcommand{\sqrpr}[1]{\left(#1\right)^2}
\newcommand{\sgnpr}[1]{\text{sign}\left(#1\right)}
\newcommand{\cjpr}[1]{\left(#1\right)^*}
\newcommand{\narop}[1]{\hspace{-0.28em}#1\hspace{-0.28em}}
\newcommand{\norm}[1]{\left|#1\right|}

\newcommand{\im}{j}

\newcommand{\Ui}{\cvec{U}_{i}}
\newcommand{\Uri}{\rvec{U}^\text{r}_{i}}
\newcommand{\Uii}{\rvec{U}^\text{i}_{i}}
\newcommand{\Uj}{\cvec{U}_{j}}
\newcommand{\Urj}{\rvec{U}^\text{r}_{j}}
\newcommand{\Uij}{\rvec{U}^\text{i}_{j}}

\newcommand{\Ilij}{\cvec{I}_{lij}}
\newcommand{\Ilji}{\cvec{I}_{lji}}
\newcommand{\Slij}{\cvec{S}_{lij}}
\newcommand{\Slji}{\cvec{S}_{lji}}
\newcommand{\Isl}{\cvec{I}^\text{s}_{lij}}
\newcommand{\Irlij}{\rvec{I}^\text{r}_{lij}}
\newcommand{\Iilij}{\rvec{I}^\text{i}_{lij}}
\newcommand{\Irlji}{\rvec{I}^\text{r}_{lji}}
\newcommand{\Iilji}{\rvec{I}^\text{i}_{lji}}
\newcommand{\Isrl}{\rvec{I}^\text{s,r}_{l}}
\newcommand{\Isil}{\rvec{I}^\text{s,i}_{l}}

\newcommand{\Zsl}{\cmat{Z}^\text{s}_{l}}
\newcommand{\Ysl}{\cmat{Y}^\text{s}_{l}}
\newcommand{\Rs}{\rmat{R}^\text{s}_{l}}
\newcommand{\Xs}{\rmat{X}^\text{s}_{l}}
\newcommand{\Yfr}{\cmat{Y}^\text{fr}_{l}}
\newcommand{\Gfr}{\rmat{G}^\text{fr}_{l}}
\newcommand{\Bfr}{\rmat{B}^\text{fr}_{l}}
\newcommand{\Yto}{\cmat{Y}^\text{to}_{l}}
\newcommand{\Gto}{\rmat{G}^\text{to}_{l}}
\newcommand{\Bto}{\rmat{B}^\text{to}_{l}}

\newcommand{\Ua}{\cscal{U}_{ia}}
\newcommand{\Ub}{\cscal{U}_{ib}}
\newcommand{\Uab}{\cscal{U}_{iab}}
\newcommand{\Uja}{\cscal{U}_{ja}}
\newcommand{\Ujb}{\cscal{U}_{jb}}
\newcommand{\Urab}{\rscal{U}^\text{r}_{iab}}
\newcommand{\Uiab}{\rscal{U}^\text{i}_{iab}}
\newcommand{\Ura}{\rscal{U}^\text{r}_{ia}}
\newcommand{\Urb}{\rscal{U}^\text{r}_{ib}}
\newcommand{\Uia}{\rscal{U}^\text{i}_{ia}}
\newcommand{\Uib}{\rscal{U}^\text{i}_{ib}}

\newcommand{\Itia}{\cscal{I}_{tijx}}
\newcommand{\Itib}{\cscal{I}_{tijy}}
\newcommand{\Itja}{\cscal{I}_{tjix}}
\newcommand{\Itjb}{\cscal{I}_{tjiy}}
\newcommand{\nt}{\rscal{n}_t}

\newcommand{\Ubagsqrab}{\rscal{U}^\text{m,sqr}_{iab}}
\newcommand{\Sd}{\cscal{S}_d}
\newcommand{\Pd}{\rscal{P}_d}
\newcommand{\Qd}{\rscal{Q}_d}
\newcommand{\Pnomd}{\rscal{P}^\text{nom}_d}
\newcommand{\Qnomd}{\rscal{Q}^\text{nom}_d}
\newcommand{\Unomd}{\rscal{U}^\text{nom}_d}
\newcommand{\gd}{\rscal{g}_d}
\newcommand{\bd}{\rscal{b}_d}

\newcommand{\Ida}{\cscal{I}_{dx}}
\newcommand{\Idb}{\cscal{I}_{dy}}
\newcommand{\Irda}{\rscal{I}^\text{r}_{dx}}
\newcommand{\Iida}{\rscal{I}^\text{i}_{dx}}
\newcommand{\Irdb}{\rscal{I}^\text{r}_{dy}}
\newcommand{\Iidb}{\rscal{I}^\text{i}_{dy}}

\newcommand{\Sg}{\cscal{S}_g}
\newcommand{\Pg}{\rscal{P}_g}
\newcommand{\Pming}{\rscal{P}^\text{min}_g}
\newcommand{\Pmaxg}{\rscal{P}^\text{max}_g}
\newcommand{\Qg}{\rscal{Q}_g}
\newcommand{\Qming}{\rscal{Q}^\text{min}_g}
\newcommand{\Qmaxg}{\rscal{Q}^\text{max}_g}
\newcommand{\Smaxg}{\rscal{S}^\text{max}_g}

\newcommand{\Iga}{\cscal{I}_{gx}}
\newcommand{\Igb}{\cscal{I}_{gy}}
\newcommand{\Irga}{\rscal{I}^\text{r}_{gx}}
\newcommand{\Iiga}{\rscal{I}^\text{i}_{gx}}
\newcommand{\Irgb}{\rscal{I}^\text{r}_{gy}}
\newcommand{\Iigb}{\rscal{I}^\text{i}_{gy}}

\newcommand{\Is}{\cvec{I}_{s}}
\newcommand{\Ss}{\cvec{S}_{s}}
\newcommand{\Irs}{\rvec{I}^\text{r}_{s}}
\newcommand{\Iis}{\rvec{I}^\text{i}_{s}}
\newcommand{\Ys}{\cmat{Y}_{s}}
\newcommand{\Gs}{\rmat{G}_{s}}
\newcommand{\Bs}{\rmat{B}_{s}}

\newcommand{\Uian}{\cscal{U}_{ian}}
\newcommand{\Urian}{\rscal{U}^\text{r}_{ian}}
\newcommand{\Uiian}{\rscal{U}^\text{i}_{ian}}
\newcommand{\Uibn}{\cscal{U}_{ibn}}
\newcommand{\Uribn}{\rscal{U}^\text{r}_{ibn}}
\newcommand{\Uiibn}{\rscal{U}^\text{i}_{ibn}}
\newcommand{\Uicn}{\cscal{U}_{icn}}
\newcommand{\Uricn}{\rscal{U}^\text{r}_{icn}}
\newcommand{\Uiicn}{\rscal{U}^\text{i}_{icn}}
\newcommand{\Unegi}{\cscal{U}^{\text{n}}_{i}}
\newcommand{\Unegri}{\rscal{U}^{\text{n,r}}_{i}}
\newcommand{\Unegii}{\rscal{U}^{\text{n,i}}_{i}}
\newcommand{\Uposi}{\cscal{U}^{\text{p}}_{i}}
\newcommand{\Uposri}{\rscal{U}^{\text{p,r}}_{i}}
\newcommand{\Uposii}{\rscal{U}^{\text{p,i}}_{i}}
\newcommand{\Uzeri}{\cscal{U}^{\text{z}}_{i}}
\newcommand{\Uzerri}{\rscal{U}^{\text{z,r}}_{i}}
\newcommand{\Uzerii}{\rscal{U}^{\text{z,i}}_{i}}
\newcommand{\Unegmaxi}{\rscal{U}^{\text{n,max}}_{i}}
\newcommand{\rot}{\alpha}
\newcommand{\vufmaxi}{\text{VUF}^\text{\,max}_{\,i}}

\newcommand{\Ilijc}{\cscal{I}_{lij,c}}
\newcommand{\Irlijc}{\rscal{I}^\text{r}_{lij,c}}
\newcommand{\Iilijc}{\rscal{I}^\text{i}_{lij,c}}
\newcommand{\Slijc}{\cscal{S}_{lij,c}}

\newcommand{\Itijc}{\cscal{I}_{tij,c}}
\newcommand{\Stijc}{\cscal{S}_{tij,c}}
\newcommand{\Irtijc}{\rscal{I}^\text{r}_{tij,c}}
\newcommand{\Iitijc}{\rscal{I}^\text{i}_{tij,c}}

\newcommand{\Sdc}{\cscal{S}_{dc}}
\newcommand{\Idc}{\cscal{I}_{dc}}
\newcommand{\Irdc}{\rscal{I}^\text{r}_{dc}}
\newcommand{\Iidc}{\rscal{I}^\text{i}_{dc}}

\newcommand{\Sgc}{\cscal{S}_{gc}}
\newcommand{\Igc}{\cscal{I}_{gc}}
\newcommand{\Irgc}{\rscal{I}^\text{r}_{gc}}
\newcommand{\Iigc}{\rscal{I}^\text{i}_{gc}}

\newcommand{\Ssc}{\cscal{S}_{sc}}
\newcommand{\Isc}{\cscal{I}_{sc}}
\newcommand{\Irsc}{\rscal{I}^\text{r}_{sc}}
\newcommand{\Iisc}{\rscal{I}^\text{i}_{sc}}

\newcommand{\Uip}{\cscal{U}_{ip}}
\newcommand{\Iliju}{\cscal{I}_{lij,u}}

\newcommand{\Idx}{\cscal{I}_{dx}}
\newcommand{\Idy}{\cscal{I}_{dy}} 
\newcommand{\Idz}{\cscal{I}_{dz}}
\newcommand{\Idw}{\cscal{I}_{dw}}

\newcommand{\Sbdx}{\cscal{S}^\text{b}_{dx}}
\newcommand{\Sbdy}{\cscal{S}^\text{b}_{dy}}
\newcommand{\Sbdz}{\cscal{S}^\text{b}_{dz}}
\newcommand{\Sbdw}{\cscal{S}^\text{b}_{dw}}
\newcommand{\Sdx}{\cscal{S}_{dx}}
\newcommand{\Sdy}{\cscal{S}_{dy}}
\newcommand{\Sdz}{\cscal{S}_{dz}}
\newcommand{\Uc}{\cscal{U}_{ic}}
\newcommand{\Un}{\cscal{U}_{in}}



\section{Introduction}

In low-voltage (LV) distribution networks, single-phase loading and generation can lead to unequal voltage magnitudes between phases. 
In networks with neutral conductors (e.g. European-style four-wire networks), the neutral current can lead to neutral-point shift, which can worsen the voltage unbalance \cite{klonari_probabilistic_2016}. 
Barbato et al. \cite{barbato_lessons_2018} report on distribution network measurements obtained through an extensive monitoring infrastructure. 
At all times, the neutral current was above 10\,\% of the peak phase current, and at several times it exceeded some of the phase currents in magnitude. 
These effects are largely overlooked in the state-of-the-art on power flow optimization in distribution networks.
Therefore, this paper compares  optimization models which are well-suited for networks where modelling the neutral wire explicitly is paramount.

\begin{figure}[tbh]
	\centering
	\includegraphics[width=3in]{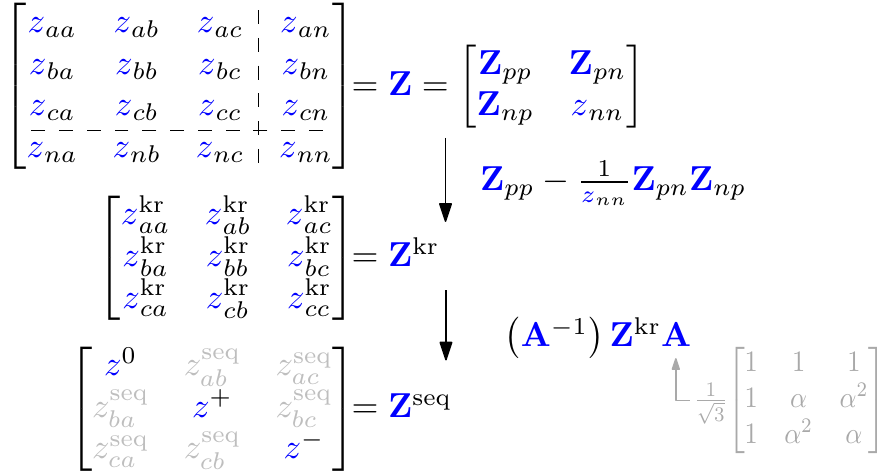}
	\caption{The approximation step applied to the series impedance. The Kron-reduction reduces the original 4$\times$4  to a 3$\times$3 matrix. The single-phase equivalent model represents the line with the complex number $\cscal{z}^+$.}
	\label{fig:brmod_approx}
\end{figure}

In the analysis of three-phase transmission networks, the network is often modeled by its positive-sequence a.k.a. `single-phase' approximation. 
In distribution networks however, lines are not transposed and line currents are unbalanced; the positive-sequence approximation is a poor substitute under these circumstances \cite{kersting_computation_2004}. 
The series impedance of a power line or cable is modeled by a square matrix of self- and mutual impedances, for each and between all conductors. 
A four-wire segment has a neutral conductor in addition to three phase conductors. 
If the neutral is perfectly grounded ($0\,\Omega$) at both ends, a `Kron-reduction' will transform the 4$\times$4 impedance matrix to an equivalent 3$\times$3 matrix. 
This simplifies the analysis, and is often done even though no such grounding exists in reality \cite{ciric_power_2003}. 
When not all nodes are actually grounded or have a nonzero grounding resistance, a Kron-reduced model will under-estimate the voltage deviations \cite{urquhart_accuracy_2016}.
The Kron-reduced $3\times3$ matrix finally is the basis for deriving the sequence components. 
The whole process is summarized in Fig. \ref{fig:brmod_approx}. There are many alternative metrics to capture unbalance levels; Girigoudar et al. discuss these in-depth, and compare them analytically \cite{girigoudar_analytical_2015}.

\subsection{State-of-the-art of Unbalanced OPF}

Ramos de Araujo et al. \cite{Araujo2013} develop an n-wire nonlinear UBOPF model in the current-voltage variable space that is capable of representing neutral wires. They use a NLP interior point method to solve the optimization problem. Dall'Anese et al. propose a bus-injection model semidefinite relaxation (BIM-SDP) of UBOPF \cite{dallanese_optimization_2012}, and Gan and Low introduced a branch-flow model semidefinite relaxation (BFM-SDP) \cite{gan_convex_2014}. 
Furthermore, there is no well-developed UBOPF benchmark library; instead, authors modify PF benchmarks ad-hoc.

Kron-reduced models assume that the return current of connected devices goes into ground at the bus. 
In four-wire networks however, the return current is injected in the neutral and - if applicable - the neutral grounding shunt. 
A nonzero current in the neutral induces voltage in other phases.
There are few studies which use four-wire network data in the context of OPF problems. 
To the best of the authors' knowledge, Su et al. \cite{su_optimal_2014} were the first to do this.  
They optimize the setpoints of PV inverters with respect to a weighted combination of network losses, generation cost, voltage deviations and unbalance. 
Recently, Gastalver-Rubio et al. \cite{gastalver-rubio_improving_2019} did a similar study. 
They minimize network losses, which is also argued to be a proxy for the unbalance. 
The UBOPF studies by Ramos de Araujo et al. \cite{Araujo2013} is arguably the most in-depth, including analysis of neutral voltage rise.

In addition to exact formulations, few authors propose convex relaxations for four-wire networks. 
Liu et al. extend previous work on BIM-SDP relaxations of UBOPF to four-wire networks \cite{liu_acopf_2019}. 
Additionally, they introduce a technique to reduce the number of variables and to recover a feasible solution when the relaxation is inexact.
The solutions are not compared with solutions obtained with an exact, non-linear formulation.
Usman et al. also proposed a four-wire extension to BIM-SDP relaxations, and further included ZIP loads \cite{usman_centralized_2019}. 
The case study examines five different test feeders, including IEEE13.
However, Usman et al. state that a wide range of solvers, including \textsc{Ipopt}, failed to solve their implementation of the exact, non-linear formulation, and therefore they cannot report results for the optimality gap.
As part of the case studies, they also examine relaxations of Kron-reduced network models, and compare the objective and tightness with the corresponding four-wire model.
In later work, they add delta-connected loads \cite{usman_bus_2020} and extensions of BFM-SDP relaxations \cite{usman_cheap_2020}.
Finally, Liu et al. \cite{liu_state_2019} also apply four-wire models in the context of state estimation, and compare rectangular and polar voltage variables.

We conclude that there is lots of ambiguous information in the literature surrounding the scalability of exact formulations of unbalanced OPF through nonlinear programming methods. 
To make sure the results are sufficiently generalizable, we must prioritize formulations which are expressive enough to support the key features used in the IEEE unbalanced power flow test cases:
\begin{itemize}
	\item transformers: lossy, multi-conductor, multi-winding;
	\item lines: multi-conductor with line-charging;
	\item loads: multi-phase, ZIP voltage dependency;
	\item shunts: multi-conductor admittance;
	\item generators: multi-phase, PQ;
	\item bounds: voltage unbalance factor, neutral shift, negative sequence voltage, phase-to-phase voltage, phase-to-neutral voltage.
\end{itemize}

\subsection{Scope and contributions}

Therefore this paper describes in detail two exact formulations to the UBOPF problem, and explores computation speed and reliability. 
Special care was taken to show how everything can be formulated as a \emph{purely quadratic} optimization problem. 
The first formulation is a quadratic one in rectangular current and voltage variables (similar to \cite{Araujo2013}); the second is the equivalent  power-voltage rectangular form. 
These formulations include many modeling features, validated against OpenDSS on several popular PF benchmarks with an accuracy of at least 1E-7.

The main contributions of this paper are 1) a methodology to convert Kron-reduced network data to four-wire data, 2) a set of numerical experiments comparing balanced, Kron-reduced and four-wire OPF models across 128 LV test feeders, in terms of constraint violation and scalability, and 3)  the observation that power balance constraints are a relaxation of Kirchhoff's current law, potentially leading to inexact solutions and numerical stability issues  four-wire OPF formulations using power balance constraints. 

\subsection{Notation and preliminaries}
The symbol for current is $I$; voltage $U$; power $S/P/Q$.
Complex variables get a blue typeface,  real variables a red one; also, vectors and matrix variables get a bold typeface.
Superscript `r' indicates the real part and `i' the imaginary part of an otherwise complex-valued symbol.

\begin{figure}[tbh]
	\centering
	\includegraphics[width=1.5in]{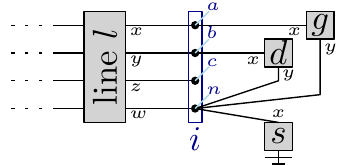}
	\caption{Illustration of a four-wire line which connects to a single-phase generator $g$, load $d$ and shunt $s$. }
	\label{fig:kcl}
\end{figure}

Figure \ref{fig:kcl} illustrates a bus $i$ that connects to a line $l$, a single-phase load $d$, a single-phase generator $g$ and a shunt $s$. We furthermore use $t$ as the transformer index. 
The buses $i$ have terminals $k \in \mathcal{P}_i \subseteq \{a,b,c,n\}$. 
The component conductors are indexed over $p\in\mathcal{U} \subseteq \{x,y,z,w\}$, e.g. for a single-phase branch $l$ we have $\mathcal{U}_l = \{x,w\}$. 
The bus voltage vector $\Ui$ stacks the terminal voltages $\Uip$, $\Ui = [\Uip]_{p \in \mathcal{P}_i}$.
The branch current vector $\Ilij$ stacks the conductor currents  $\Ilij = [\Iliju]_{u \in \mathcal{U}_l}$. 
We therefore define an mapping for each component, from its conductors to (a subset of) the bus terminals. For instance for the branch with $\mathcal{U}_l= \{x,w \}$ we have the from-side conductor-terminal mapping $\mathcal{X}_{lij} = \{(x,b),(w,n)\} \subset \mathcal{U}_l \times \mathcal{P}_i$.
Note that this data model allows for arbitrary connections without relying on the order of the conductors.

\begin{figure}[tbh]
	\centering
	\subfloat[wye-connected]{\includegraphics[width=1.6in]{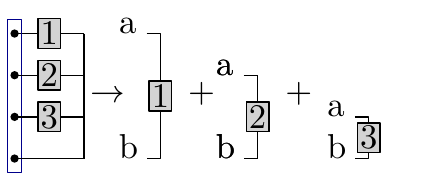}%
	\label{fig:config_wye}}
	\hspace{0.2em}
	\subfloat[delta-connected]{\includegraphics[width=1.6in]{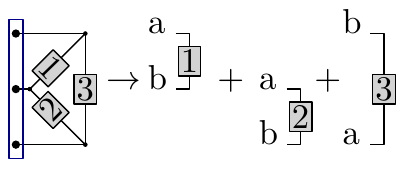}%
	\label{fig:config_delta}}
	\caption{Three-phase components can be equivalently composed of three single-phase components.}
	\label{fig:config}
\end{figure}

Where-ever more concise, we will make abstraction of delta or wye connections for devices. 
Figure \ref{fig:config} shows that both can be interpreted as equivalent compositions of single-phase components.
For transformers, which are two-port components, this can be done for each port separately, leading to Yy, Dd, Yd and Dy transformers.

\section{Current-voltage Formulation}\label{sect:form}

The derivation uses primarily complex current (I) and voltage (V) variables, both split into rectangular (R) coordinates; hence, it is referred to as the IVR formulation.
\begin{figure}[tbh]
	\centering
	\subfloat[multi-conductor line]{\includegraphics[width=1.6in]{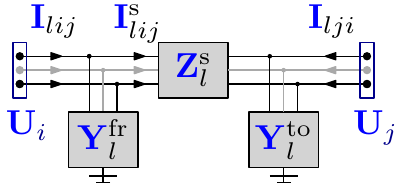}%
	\label{fig:branch}}
	\hspace{0.5em}
	\subfloat[single-phase transformer]{\includegraphics[width=1.7in]{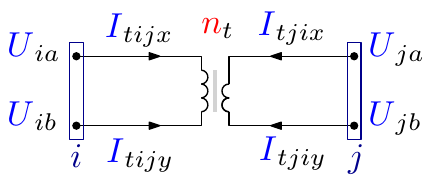}%
	\label{fig:transformer}}
	\caption{Two-port (excluding ground) components connect two buses $i$ and $j$.}
	\label{fig:twoport}
\end{figure}

\subsection{Multi-conductor line}
A multi-conductor line $l$ consists of $n$ conductors, which connect $n$ terminals at bus $i$ to $n$ terminals at bus $j$.
A $\pi$-model of the line is parameterized by three complex matrices, $\Zsl$, $\Yfr$ and $\Yto$.
Figure \ref{fig:branch} illustrates the $\pi$-model, governed by
\begin{subequations}\begin{align}
	\Ui-\Uj &= \Zsl\Isl \label{eqc:line_ohm},\\
	\Ilij = \Yfr\Ui + \Isl &, \quad
	\Ilji = \Yto\Uj - \Isl \label{eqc:line_cto}.
\end{align}\end{subequations}

\subsection{Single-phase transformer}
Figure \ref{fig:transformer} illustrates an ideal, two-winding, single-phase transformer $t$, with turns ratio $\nt$. Note that $\nt$ is a real and therefore does not model any phase shifting. 
We observe the conductor-terminal maps  $\mathcal{X}_{lij} = \{(x,a),(y,b)\}$, $\mathcal{X}_{lji} = \{(x,a),(y,b)\}$.
The transformer's feasible set is therefore,
\begin{subequations}\begin{align}
	\Itia + \Itib &= 0,\\
	\Itja + \Itjb &= 0,\\
	\nt\Itia + \Itja &= 0,\\
	\Ua-\Ub &= \nt\left(\Uja-\Ujb\right).
\end{align}\end{subequations}
A multi-winding, multi-phase, lossy transformer can be equivalently represented by
a composition of lines, shunts and ideal two-winding transformers \cite{Claeys2020b}.
Combining these  elementary building blocks, this formulation can represent all transformers covered in \cite{Claeys2020b}.

\begin{figure}[tbh]
	\centering
	\subfloat[load]{\includegraphics[width=0.8in]{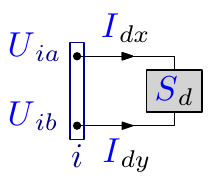}%
	\label{fig:load}}
	\hspace{1em}
	\subfloat[generator]{\includegraphics[width=0.8in]{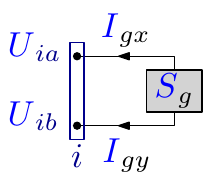}%
	\label{fig:gen}}
	\hspace{1em}
	\subfloat[shunt]{\includegraphics[width=0.6in]{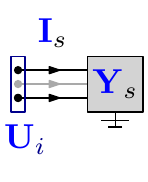}%
	\label{fig:shunt}}
	\caption{One-port components connect to a single bus.}
	\label{fig:oneport}
\end{figure}

\subsection{Load}
\subsubsection{Three-phase}

A three-phase, wye-connected load connects to 
bus $i$ with the conductor-terminal map $\{(x,a),(y,b),(z,c),(w,n)\}$. The load current and power variables therefore satisfy,
\begin{subequations}\label{eqc:load_3ph}\begin{align}
 0 &= \Idx + \Idy + \Idz + \Idw,\\
 \Sdx &= (\Ua-\Un)\cjpr{\Idx},\\
 \Sdy &= (\Ub-\Un)\cjpr{\Idy},\\
 \Sdz &= (\Uc-\Un)\cjpr{\Idz},
\end{align}\end{subequations}
where e.g. $\Sdx$ is the power drawn between terminals $a$ and $n$.
The power flow  for conductor $x$ is defined as $\Sbdx=\Ua\cjpr{\Idx}$, and similarly for the others.
To define voltage-sensitive loads, we consider three-phase loads as the composition of single-phase loads. 

\subsubsection{Single-phase}

Figure \ref{fig:load} illustrates a single-phase load $d$, which draws an amount of power $\Sd$ between the terminals $a$ and $b$ at bus $i$, with map  $\mathcal{X}_{di} = \{(x,a),(y,b)\}$,
\begin{subequations}\begin{align}
	\Ida + \Idb &= 0 \label{eqc:load_cbal},\\
	\Sd &= \left(\Ua-\Ub\right)\cjpr{\Ida} \label{eqc:load_sdef}.
\end{align}\end{subequations}
Loads are further characterized by how the drawn power depends on the voltage magnitude applied to them, $\Sd = f(\norm{\Ua-\Ub})$. Below, several common models are discussed. We note that ZIP loads can be represented as a composition of a constant impedance, constant current and a constant power load.

\paragraph{Constant power} The power is constant and independent of the applied voltage.
\begin{align}
	\Pd &= \Pnomd, & \Qd &= \Qnomd. \label{eq_fix_power_load}
\end{align}

\paragraph{Constant impedance} The power is proportional to the square of the voltage magnitude,
\begin{align}
	\Pd &= \frac{\Pnomd}{\sqrpr{\Unomd}}\norm{\Uab}^2, & \Qd &= \frac{\Qnomd}{\sqrpr{\Unomd}}\norm{\Uab}^2. \label{eq_const_S_power}
\end{align}
where we introduced the shorthand $\Uab=\Ua-\Ub$. This seems like it is quadratic, but can be represented linearly in the variables $\Uab$ and $\Ida$. If $\Uab\neq 0$,
\begin{subequations}\label{eqc:load_z}\begin{align}
 \Leftrightarrow\; &  \Ida=\left(\gd+\im\bd\right)\Uab, \label{eq_const_S_current}\\
									  \hspace{1em} \gd = \hspace{0.7em}\Pnomd &/ \sqrpr{\Unomd},
									  \hspace{1em} \bd = -\Qnomd / \sqrpr{\Unomd}. 
\end{align}\end{subequations}

\paragraph{Constant current} The power is proportional to the voltage magnitude,
\begin{align}
	\Pd &= \frac{\Pnomd}{\Unomd}\norm{\Uab}, & \Qd &= \frac{\Qnomd}{\Unomd}\norm{\Uab}.
\end{align}
This can be represented using an auxiliary variable $\Ubagsqrab$ defined as, 
\begin{align}
	\Ubagsqrab = \left|\Ua\narop{-}\Ub\right|^2 = \left(\Ura\narop{-}\Urb\right)^2+\left(\Uia\narop{-}\Uib\right)^2, \label{eq_const_current_aux}
\end{align}
and two quadratic equalities and two linear inequalities, 
\begin{subequations}\begin{align}
	\sqrpr{\Pd} &= \sqrpr{\frac{\Pnomd}{\Unomd}} \Ubagsqrab, & \sgnpr{\Pnomd}\Pd \geq 0, \label{eq_const_current_P}\\
	\sqrpr{\Qd} &= \sqrpr{\frac{\Qnomd}{\Unomd}} \Ubagsqrab, & \sgnpr{\Qnomd}\Qd \geq 0. \label{eq_const_current_Q}
\end{align}\end{subequations}


\subsection{Single-phase generator}
Figure \ref{fig:gen} shows a single-phase generator $g$, with map  $\mathcal{X}_{gi} = \{(x,a),(y,b)\}$,
\begin{subequations}\begin{align}
	\Iga + \Igb &= 0, \label{eqc:gen_cbal}\\
	\Sg &= \left(\Ua-\Ub\right)\cjpr{\Iga}. \label{eqc:gen_sdef}
\end{align}\end{subequations}
The power $\Sg=\Pg+\im\Qg$ it injects, is bounded by
\begin{subequations}\label{eq:gen_bounds}\begin{align}
	\Pming \leq \Pg \leq \Pmaxg, \quad &
	\Qming \leq \Qg \leq \Qmaxg, \label{eq_gen_bounds_PQ}\\
	\sqrpr{\Pg} + \sqrpr{\Qg} &\leq \sqrpr{\Smaxg}. \label{eq_gen_bounds_S}
\end{align}\end{subequations}

\subsection{Shunt}
Figure \ref{fig:shunt} shows a multi-conductor shunt $s$, defined by,
\begin{align}
	\Is &= \Ys\Ui .\label{eqc:shunt}
\end{align}
Lossy grounding of the neutral conductor can be modelled as a single-conductor shunt.

\subsection{Kirchhoff's current law at buses}
Let $\mathcal{L}_{ik}$ denote the set of all line conductors $lijp$ connected to terminal $(i,k)$ at bus $i$; the same goes for $\mathcal{T}_{ik}$ (transformers), $\mathcal{D}_{ik}$ (loads), $\mathcal{S}_{ik}$ (shunts) and $\mathcal{G}_{ik}$ (generators). The expression for Kirchhoff's current law, requiring that the sum of all incoming currents equals the sum of all outgoing currents at each terminal in the network, at each terminal $(i,k)$,
\begin{align}
&\sum_{gp\in\mathcal{G}_{ik}} \!\!\Igc = \!\!\!\sum_{lijp\in\mathcal{N}_{ik}} \!\! \Ilijc +\!\!\sum_{tijp\in\mathcal{T}_{ik}} \!\!\Itijc +\!\! \sum_{dp\in\mathcal{D}_{ik}}\!\! \Idc + \!\! \sum_{sp\in\mathcal{S}_{ik}} \!\! \Isc
.\label{eqc:kcl}\end{align}
For the other terminals at the bus, we define similar sets and expressions. 
Figure \ref{fig:kcl} illustrates this with an example, where a bus connects to a line $l$, a single-phase load $d$, a single-phase generator $g$ and a shunt $s$. Applying (\ref{eqc:kcl}) at each of the four terminals,
\begin{subequations}\begin{align}
	\cscal{I}_{gx} 	&= \cscal{I}_{lijx}, & 0 							&= \cscal{I}_{lijy} + \cscal{I}_{dx},\\
	0 							&= \cscal{I}_{lijz}, & \cscal{I}_{gy} &= \cscal{I}_{lijw} + \cscal{I}_{dy} + \cscal{I}_{sx}.
\end{align}\end{subequations}

\subsection{Technical envelopes}

\paragraph{Current ratings}
The magnitude of each of the current variables $\cscal{I}=\rscal{I}^\text{r}+\im\rscal{I}^\text{i}$ introduced before,
can be enforced with a second-order-cone (SOC) constraint,
\begin{align}
	(\rscal{I}^\text{r})^2 + (\rscal{I}^\text{i})^2 \leq (\rscal{I}^\text{max})^2. \label{eq_current_bound_current}
\end{align}

\paragraph{Voltage magnitude bounds}
Some voltage magnitude bounds apply to the difference of two voltage phasors (phase-to-phase $\rscal{U}^\text{min}_{iab}, \rscal{U}^\text{max}_{iab}$ and phase-to-neutral bounds $\rscal{U}^\text{min}_{ian}, \rscal{U}^\text{max}_{ian}$), whilst others apply directly to a single voltage phasor (neutral shift bound, $\rscal{U}^\text{min}_{in}, \rscal{U}^\text{max}_{in}$). Both can be expressed as quadratic inequalities,
\begin{align}
	&(\rscal{U}^\text{min}_{ia})^2 &&\leq &&(\Ura)^2 &&+ &&(\Uia)^2 &&\leq &&(\rscal{U}^\text{max}_{ia})^2, \label{eq_volt_magn_bound}\\
	&(\rscal{U}^\text{min}_{iab})^2 &&\leq &&(\Urab)^2 &&+ &&(\Uiab)^2 &&\leq &&(\rscal{U}^\text{max}_{iab})^2, \label{eq_volt_diff_magn_bound}
\end{align}
where the upper bound is representable as a SOC constraint, whilst the lower bound is non-convex.

\paragraph{Unbalance bounds}
Take a bus $i$, with phase terminals $a$, $b$ and $c$, and neutral $n$. We assume $b$ lags $a$, $c$ lags $b$. The positive sequence $\Uposi$, negative sequence $\Unegi$ and zero sequence $\Uzeri$ phasors are then defined,
\newcommand{\coeffspacing}{-0.8em}
\begin{subequations}\begin{align}
&\Uzeri = \Uzerri + \im \Uzerii= &\tfrac{1}{3} &&\hspace{\coeffspacing}\Uian &&+ &&\tfrac{1}{3} &&\hspace{\coeffspacing}\Uibn &&+ &&\tfrac{1}{3} &&\hspace{\coeffspacing}\Uicn\label{eqc:uzer},\\
&\Uposi = \Uposri + \im\Uposii= &\tfrac{1}{3} &&\hspace{\coeffspacing}\Uian &&+ &&\tfrac{\rot}{3} &&\hspace{\coeffspacing}\Uibn &&+ &&\tfrac{\rot^2}{3} &&\hspace{\coeffspacing}\Uicn\label{eqc:upos},\\
&\Unegi = \Unegri + \im\Unegii = &\tfrac{1}{3} &&\hspace{\coeffspacing}\Uian &&+ &&\tfrac{\rot^2}{3} &&\hspace{\coeffspacing}\Uibn &&+ &&\tfrac{\rot}{3} &&\hspace{\coeffspacing}\Uicn\label{eqc:uneg},
\end{align}\end{subequations}
where $\rot=\angle 120^\circ = -0.5+\im\tfrac{\sqrt{3}}{2}$.
The voltage unbalance factor, $\text{VUF}=|\Unegi|/|\Uposi|$, is a common measure for the amount of `unbalance'.
Some grid codes specify $\text{VUF}_i\leq\vufmaxi$,
\begin{align}
(\Unegri)^2 + (\Unegii)^2  &\leq (\vufmaxi)^2\left((\Uposri)^2 + (\Uposii)^2\right), \label{eq_vuf_bound}
\end{align}
whilst others limit the negative sequence, $|\Unegi|\leq\Unegmaxi$,
\begin{align}
(\Unegri)^2 + (\Unegii)^2  &\leq (\Unegmaxi)^2. \label{eq_neg_seq_bound}
\end{align}

\paragraph{Power bounds}
Where desired,  bounds on the product of any $\cscal{U}$ and $\cscal{I}$ can be imposed quadratically by introducing auxiliary variables $\rscal{P}$ and $\rscal{Q}$, and imposing the constraints
\begin{align}
	\rscal{P} &= \rscal{U}^\text{r}\rscal{I}^\text{r} + \rscal{U}^\text{i}\rscal{I}^\text{i}, &
	\rscal{Q} &= -\rscal{U}^\text{r}\rscal{I}^\text{i} + \rscal{U}^\text{i}\rscal{I}^\text{r}.
\end{align}
Analogous to (\ref{eq:gen_bounds}), bounds $\rscal{P}^\text{min}$, $\rscal{P}^\text{max}$, $\rscal{Q}^\text{min}$, $\rscal{Q}^\text{max}$ and $\rscal{S}^\text{max}$ can then be imposed.

\section{Power-Voltage Formulation}\label{sect:form_SU}
In this section we present derivations that primarily use \emph{complex power } and voltage variables, again both split into rectangular (R) coordinates; hence, it is referred to as the ACR form.
Because IVR and ACR share rectangular voltage variables, all expressions that just depend on voltage can be used as-is. For brevity, we don't repeat such constraints in this section, but list the shared constraints later when defining the feasible sets (\S \ref{sect:formulation_comparison}). 


\subsection{Kirchhoff's current law at buses}
We reformulate Kirchhoff's current law, $\sum \cscal{I}=0$, by multiplying its conjugate on both sides by the nodal voltage $\cscal{U}$, obtaining $\sum \cscal{U}\cscal{I}^*=0$. Introducing variables $\cscal{S}=\cscal{U}(\cscal{I})^*$ and substituting these, then leads to $\sum \cscal{S}=0$, the power balance form. The `flow variable' interfacing the components, now is power instead of current:

\begin{align}
&\sum_{gp\in\mathcal{G}_{ik}} \!\!\Sgc = \!\!\!\sum_{lijp\in\mathcal{N}_{ik}} \!\! \Slijc +\!\!\sum_{tijp\in\mathcal{T}_{ik}} \!\!\Stijc +\!\! \sum_{dp\in\mathcal{D}_{ik}}\!\! \Sdc + \!\! \sum_{sp\in\mathcal{S}_{ik}} \!\! \Ssc
.\label{eqc:kcl_power}\end{align}

\subsection{Multi-conductor line}
Lines and shunts allow a compact representation in terms of power flow variables. The from and to-side line flow variables, $\Slij$ and $\Slji$, are linked to the voltage,
\begin{align}
		\Slij &= \Ui\odot\left(\Yfr\Ui+\Ysl(\Ui-\Uj)\right)^*,\label{eq_SU_power_from_def} \\
		\Slji &= \Uj\odot\left(\Yto\Uj+\Ysl(\Uj-\Ui)\right)^*, \label{eq_SU_power_to_def}
\end{align}
where $\Ysl=(\Zsl)^{-1}$, and $\odot$ denotes the element-wise product. 

\subsection{Shunt}
The shunt power flow variable, $\Ss$, is linked to the voltage by
\begin{align}
		\Ss = \Ui\odot(\Ys\Ui)^*. \label{eq_shunt_power}
\end{align}
\subsection{Loads, generators and transformers}
We redefine the load model with power flow variables by reformulating (\ref{eqc:load_3ph}) with $\Sbdx=\Ua\cjpr{\Idx}$ etc.,
\begin{subequations}\label{eqc:load_3ph_p}\begin{align}
0 &=
\begin{aligned}
&\Ub\Uc\Un\Sbdx + \Ua\Uc\Un\Sbdy\\
+\;&\Ua\Ub\Un\Sbdz + \Ua\Ub\Uc\Sbdw, 
\end{aligned}
\label{eqc:load_pflowkcl}\\
 \Sdx\Ua &= \Sbdx(\Ua-\Un),\\
 \Sdy\Ua &= \Sbdy(\Ua-\Un),\\
 \Sdz\Ua &= \Sbdz(\Ua-\Un).
\end{align}\end{subequations}
Note that (\ref{eqc:load_pflowkcl}) requires the introduction of several auxiliary variables to reformulate it as a set of quadratic constraints (i.e. repeated substitution $z\leftarrow xy$). 
This approach, despite including a fourth degree polynomial, is popular for Kron-reduced forms, because when $\Un=0$ and terminal $n$ is grounded, (\ref{eqc:load_3ph_p}) reduces to $\Sbdx=\Sdx$, $\Sbdy=\Sdy$ and $\Sbdz=\Sdz$.

Generators and transformers are similarly reformulated in terms of power flow variables.

\subsection{Technical envelopes}

\paragraph{Current ratings}
The magnitude of each of the current variables $\cscal{I}=\rscal{I}^\text{r}+\im\rscal{I}^\text{i}$ introduced before,
can be constrained with a second-order-cone (SOC) constraint in the SU variable space,
\begin{align}
	(\rscal{P})^2  + (\rscal{Q})^2 \leq (\rscal{I}^\text{max})^2 \left( (\rscal{U}^{\text{r}})^2 + (\rscal{U}^{\text{i}})^2 \right). \label{eq_current_bound_current_lifted}
\end{align}

\section{Assembling four-wire OPF data sets}

To obtain a realistic amount of unbalance, we use the dataset released by Rigoni et al., i.e. a set of 128 LV, four-wire distribution networks, representative of feeders in the UK \cite{rigoni_representative_2016}.
However, only the positive and zero sequence components of the  lines are included, whilst for the upcoming numerical experiments,  full 4$\times$4 impedance data is paramount.
We now present a methodology to obtain the missing 4$\times$4 data, and to derive OPF problems with interesting dispatchable resources and envelopes.

\subsection{Obtaining four-wire data}
The original dataset only includes positive and zero sequence components.
Figure \ref{fig:brmod_approx} shows how these are calculated based on the 4$\times$4 impedance matrices.

Carson's equations provide a method to calculate the 4$\times$4 matrices, based on the spacing and properties of the conductors and the ground \cite{kersting_computation_2004}.
A linecode specifies the series impedance of a line per unit of length.
The following procedure assigns to each line an alternative linecode, which contains four-wire impedance data:
\begin{enumerate}
	\item Compile a library of candidate linecodes $c\in\mathcal{C}$, of which the 4$\times$4 impedance matrices are known.
	\item For each of these, calculate $\cscal{z}^+_c=\rscal{r}^+_c+\im\rscal{x}^+_c$.
	\item For each linecode $o$ in the original dataset, replace it with the linecode that is closest to it in ($\rscal{r}^+$,$\rscal{x}^+$) space,
	\begin{align}
			o \leftarrow \underset{c\in\mathcal{C}}{\operatorname{argmin}} \begin{vmatrix}r^+_o-r^+_c\\x^+_o-x^+_c\end{vmatrix}_2.
	\end{align}
\end{enumerate}

\begin{figure}[tbh]
	\centering
	\includegraphics[width=3.5in]{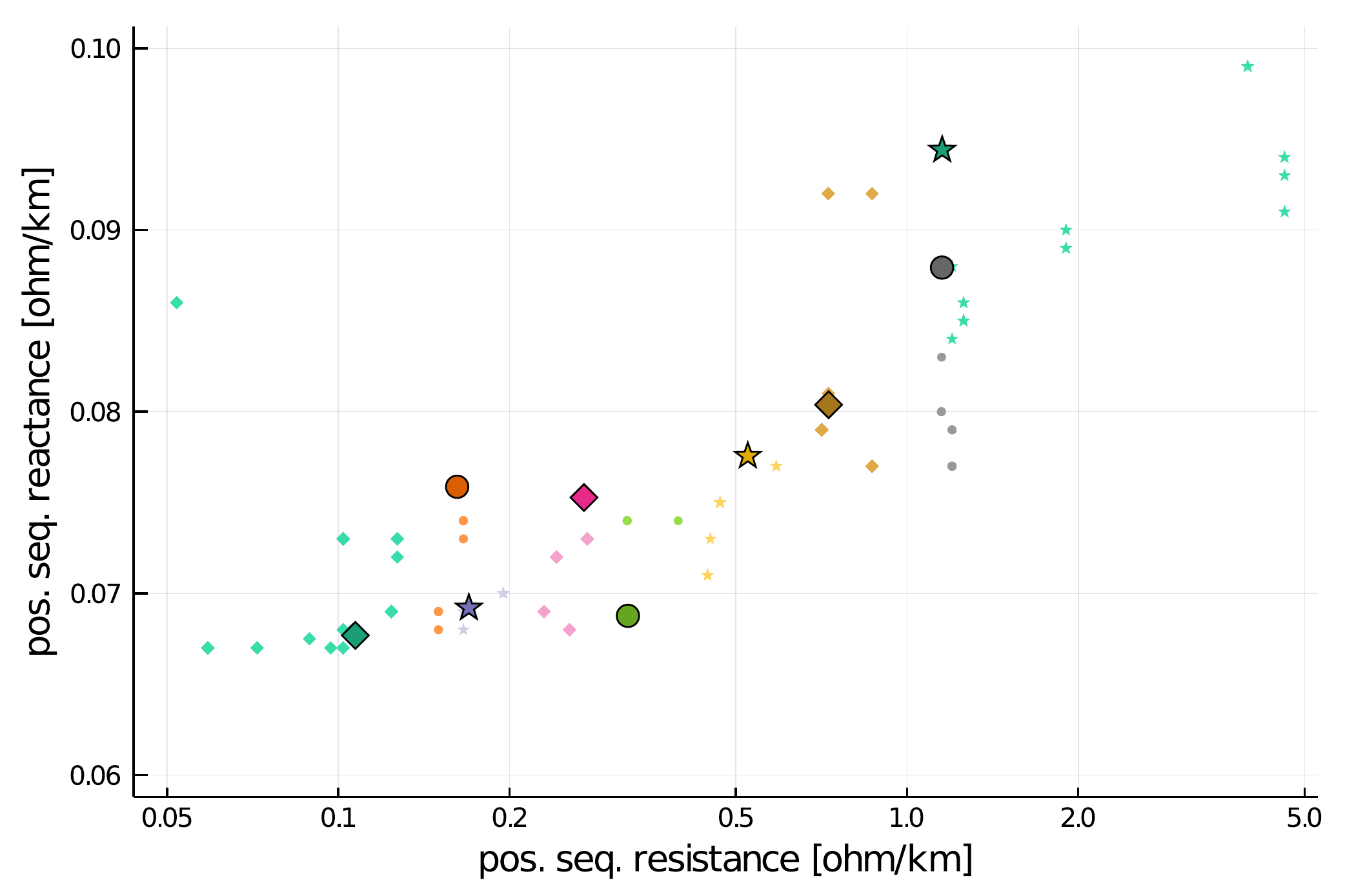}\\[-0.6em]
	\caption{Each original linecode (small markers) is replaced by a candidate linecode (big markers). Note that the x-axis has a logarithmic scale, which makes the clustering seem a bit counter-intuitive in this figure. One original linecode outlier, ($r^+_o=0.162$, $x^+_o=0.050$), is not displayed in this plot.}
	\label{fig:cluster}
\end{figure}

Urquhart et al. studied in depth the assumptions and approximations typically made in simulating low-voltage distribution networks.
As part of this study, the full impedance data for several underground linecodes common in the UK, were calculated \cite{urquhart_series_2015}.
Linecodes corresponding to conductors with cross section areas of $16$ (x2), $25$, $35$, $70$, $95$, $120$, $185$ and $300\,mm^2$, are used to populate the library of candidate linecodes.
Figure \ref{fig:cluster} shows graphically how each linecode in the original dataset was mapped to one of these 9 candidate linecodes.

Once assigned, the candidate linecodes are used to re-derive the Kron-reduced and single-phase equivalent models, as shown in Figure \ref{fig:brmod}. The neutral in the four-wire models, is grounded only at the substation (lossless). Also, the original data contains sections with many short lines  series; this is a consequence of tracing the network topology from  GIS data. When these short lines have zero shunt admittance, and no bounds are imposed on the lines or buses connecting them, the network can be simplified. Segments of short lines in series can be replaced by a single long line, leading to fewer buses but identical results. The network sizes reported in this paper are the ones after aforementioned reduction.

\subsection{Setting up interesting OPF cases}
The data set released by Rigoni et al. \cite{rigoni_representative_2016} contains 100 time series for various devices (households, electric vehicles, etc.), with a horizon of 24 h and a resolution of 5 min. We apply the load and photovoltaic time series with a power factor of 0.95. This case study will compare dispatch results for a single snapshot. We choose timestep $t=207$ such that the phase-to-neutral voltage magnitude is within the upper bound. The upstream grid is modelled as a voltage source: the voltage phasors are fixed at the reference bus, and an unbounded generator models the import/export.

Next, for each network separately, we add a set of dispatchable, distributed generators (DGs) which are allocated and sized based on the following procedure.
\begin{itemize}
	\item Add a single-phase DG at every 4th load. Each DG has a generation cost which is ten times smaller than the generation cost of the reference bus generator, modelling the grid import. Intuitively, this will lead the DG output to increase until some operational bound becomes active.
	\item Each DG $g$ has the same capacity $\Pmaxg$ (and furthermore $\Pming=\Qming=0$, $\Qmaxg=\Pmaxg$), which is sized such that under the balanced decision model, the DGs can congest the network.
\end{itemize}

The objective is to minimize the total DG dispatch cost. 
The main operational constraints are the generator bounds, and upper and lower bounds on the phase-to-neutral voltage magnitude at buses which have loads or generators connected.

\section{Implementation and Validation}
The formulations presented  focus on simplicity of the derivation rather than
efficiency in terms of number of optimization variables. 
Therefore this subsection
discusses implementation details that were left out before.
\paragraph{Complex to real}
Most optimization modelling languages do not support complex parameters and variables,
and require the user to enter the problem in real variables. 
Most expressions are straight-forward generalisations of those presented in \cite{geth_current-voltage_2020}. 
\paragraph{Perfect grounding}\label{sect:pfg}
When a terminal $(i,a)$ is perfectly grounded, the associated voltage variables can be
replaced by constants, $\Uia\leftarrow0$, $\Ura\leftarrow 0$. Instead of introducing
a variable for the grounding current, which would act as a slack variable in KCL,
one can drop the KCL constraint for $(i,a)$ and recover the grounding current in post-processing.

\paragraph{Unnecessary conductor or terminal variables}
We do not assume the network is 4-wire throughout and/or all buses have four 4 terminals. Instead of eliminating those (free) variables, e.g. by fixing them to 0 numerically, we chose not to define the variables to begin with. 
Because of this, extensions to n-conductor forms and/or multiple neutrals should be straight-forward and numerically efficient.

\paragraph{Efficient three-phase composition}
A three-phase, wye-connected load only requires 4 current variables.
However, when it is composed of three equivalent single-phase loads, 6 current variables are used.
In the implementation, we did the proper substitutions to eliminate these extra
2 current variables, and more generally as well for other multi-phase components.

\paragraph{Solver}
All OPF instances were solved with the primal-dual interior point solver \textsc{Ipopt} \cite{wachter_implementation_2006} (version 3.13.2, with MUMPS) using the automatic differentiation as implemented in JuMP/MOI \cite{Legat2022}. 

\paragraph{Initialization} Proper initialization of the voltage variables was crucial to achieve reliable convergence with \textsc{Ipopt}. The voltage variables were initialized at the value they would take under no-load conditions (and ignoring any shunts in the system), i.e. we assume a balanced phasor at the reference bus and propagate that across transformers to include the vector group angle shift, while keeping track of the terminal-conductor maps.  

\paragraph{Validation}
The proposed formulation can also be used to solve PF problems as feasibility problems.
This approach was used to validate the optimization models against the simulation models included in OpenDSS.
Table \ref{tab:validation} shows that across several well-known benchmarks for
unbalanced PF, the largest relative deviation for any
voltage phasor is 1.2E-8.
These results are obtained with \textsc{Ipopt} (version 3.13.2, tolerance set to 1E-10),
and compared against OpenDSS with tolerance set to 1E-10.


\begin{table}[h!]
	\caption{The max. relative error on the voltage  does not exceed 1.2E-8.}
	\label{tab:validation}
	\begin{center}
	\begin{tabular}{lllll}\toprule
		                    & IEEE13 & IEEE34 & IEEE123$^*$ & LVTestCase (t=1000)\\\midrule
		$\delta^\text{max}$ & 2.8E-8 & 7.7E-8 & 1.2E-8 & 3.4E-8\\
		\bottomrule
	\end{tabular}\end{center}
	\hspace{2em}$^*$: Bus 610 is compared phase-to-phase.
\end{table}



\section{Numerical Experiments}

This section discusses the constraint violation and solve time of approximate balanced and Kron-reduced models, compared to explicit four-wire models.

\subsection{Formulation variants and calculation time}\label{sect:formulation_comparison}
The results are obtained with explicit variables for the total and series line current. However, the total current variables can be replaced by affine transformations of the series current to make the formulation more efficient. For example, for the four-wire network model, this further reduces the solve time on average by 40\% (note that  line current / power are unbounded).  
Bounds are often useful to keep thee solver from diverging; for these cases, voltage magnitude bounds alone were sufficient to do so. The feasible sets are defined by \S\ref{sect:form} and \S\ref{sect:form_SU} for the IVR and ACR formulations.


\subsection{Methodology for comparing UBOPF forms}
\begin{figure}[tb]
	\centering
	\subfloat[four-wire]{\includegraphics[width=0.8in]{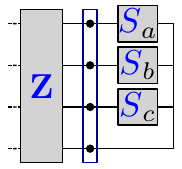}%
	\label{fig:brmod_fw}}
	\hspace{1em}
	\subfloat[Kron-reduced]{\includegraphics[width=0.8in]{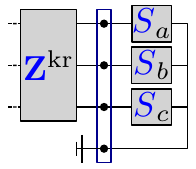}%
	\label{fig:brmod_kr}}
	\hspace{1em}
	\subfloat[single-phase equivalent]{\includegraphics[width=1.2in]{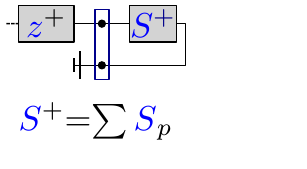}%
	\label{fig:brmod_bal}}
	\caption{Illustration of different network models, from more to less detailed.}
	\label{fig:brmod}
\end{figure}

The goal of the numerical experiments is to compare optimization results obtained with three different models: `four-wire', `Kron-reduced' and  `single-phase equivalent'. 
Figure \ref{fig:brmod} illustrates how each of these can be represented by the proposed formulation.

\begin{figure}[tb]
	\centering
	\includegraphics[width=2.2in]{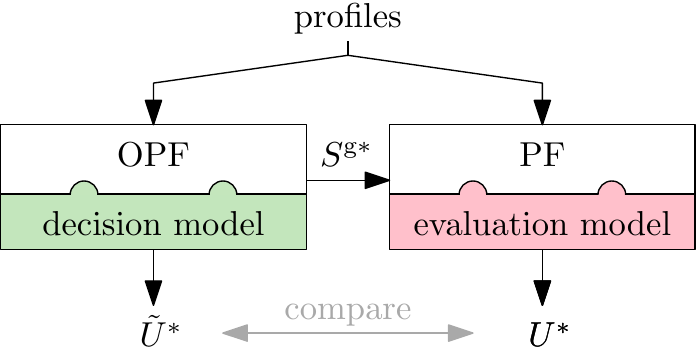}
	\caption{The dispatch results obtained with a less detailed OPF model (Kron-reduced or single-phase equivalent), are evaluated with a four-wire PF model. The voltages can then be compared, quantifying the approximation error.}
	\label{fig:evaluation}
\end{figure}

Figure \ref{fig:evaluation} shows how the accuracy of the simplified models (Kron-reduced and single-phase equivalent) is evaluated.
Once the dispatch results are obtained, the `true' voltage profile is then obtained by solving a four-wire PF \emph{simulation} problem. 
The difference between these models typically increases as the networks are more unbalanced.

\subsection{Current-voltage formulation}
\subsubsection{Constraint violation}

\begin{figure}[tbh]
	\centering
	\def\arraystretch{0.0}
	\subfloat[balanced]{\begin{tabular}[b]{c}\includegraphics[width=1.5in,trim=0 0.2in 0 0.1in,clip]{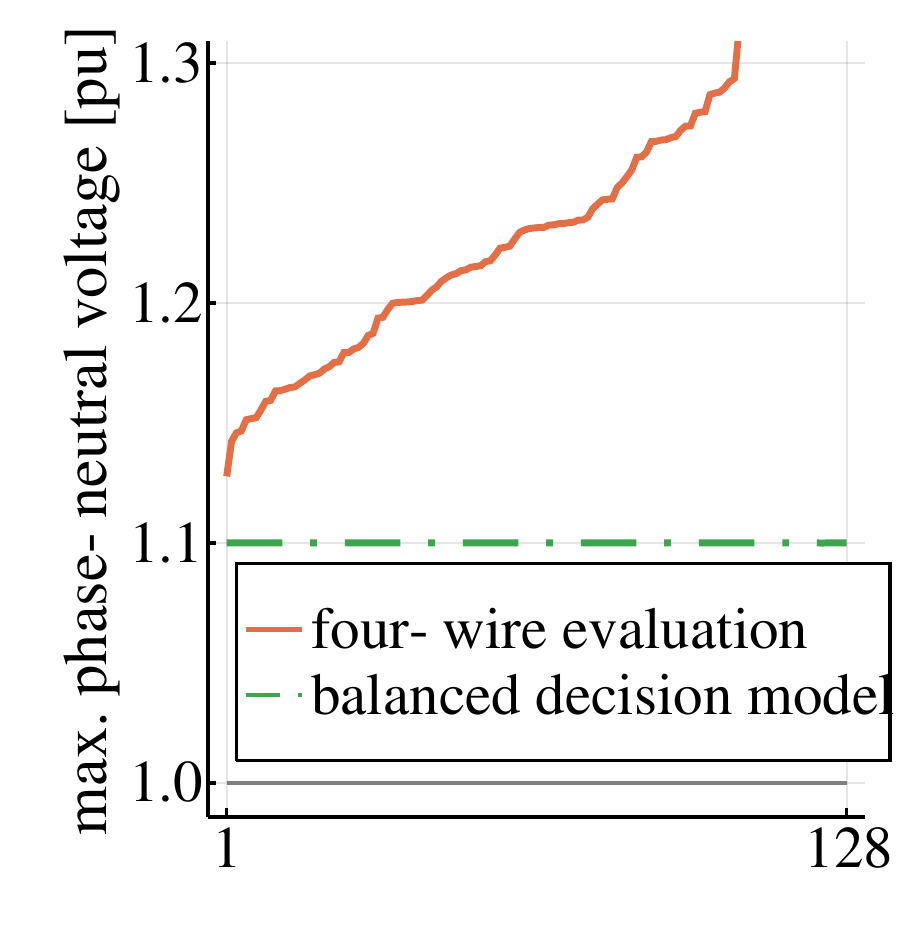}\\
	\hspace{2em}\footnotesize{network instance}\end{tabular}\label{fig:pnmax_10_bal}}
	\hspace{0em}
	\subfloat[Kron-reduced]{\begin{tabular}[b]{c}\includegraphics[width=1.5in,trim=0 0.2in 0 0.1in,clip]{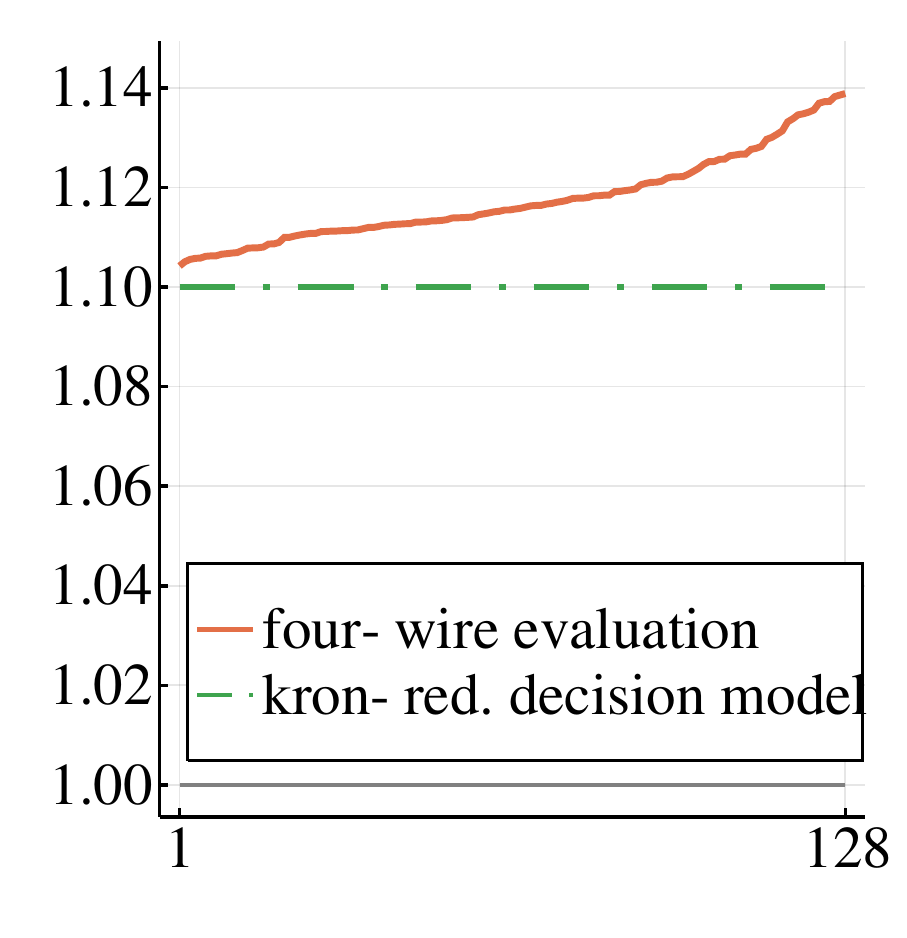}\\
	\hspace{2em}\footnotesize{network instance}\end{tabular}\label{fig:pnmax_10_kron}}
	\caption{Evaluation of the highest constraint violation (phase-to-neutral limited to 1.10 pu) caused by approximate network models, for all 128 instances (ordered by violation amount).}
	\label{fig:pnmax_10}
\end{figure}

Following the methodology in Figure \ref{fig:evaluation}, the approximate balanced and Kron-reduced models are used to optimize the generator setpoints. These setpoints are then evaluated by solving the four-wire PF. 
Figure \ref{fig:pnmax_10} shows the highest phase-to-neutral voltage magnitude that this would cause for each of the 128 network instances, for an imposed upper bound of 1.1 pu. 
Even using the Kron-reduced model, leads to a constraint violation of 1-4\%.

\begin{figure}[tbh]
	\centering
	\def\arraystretch{0.0}
	\subfloat[balanced]{\begin{tabular}[b]{c}\includegraphics[width=1.5in,trim=0 0.2in 0 0.1in,clip]{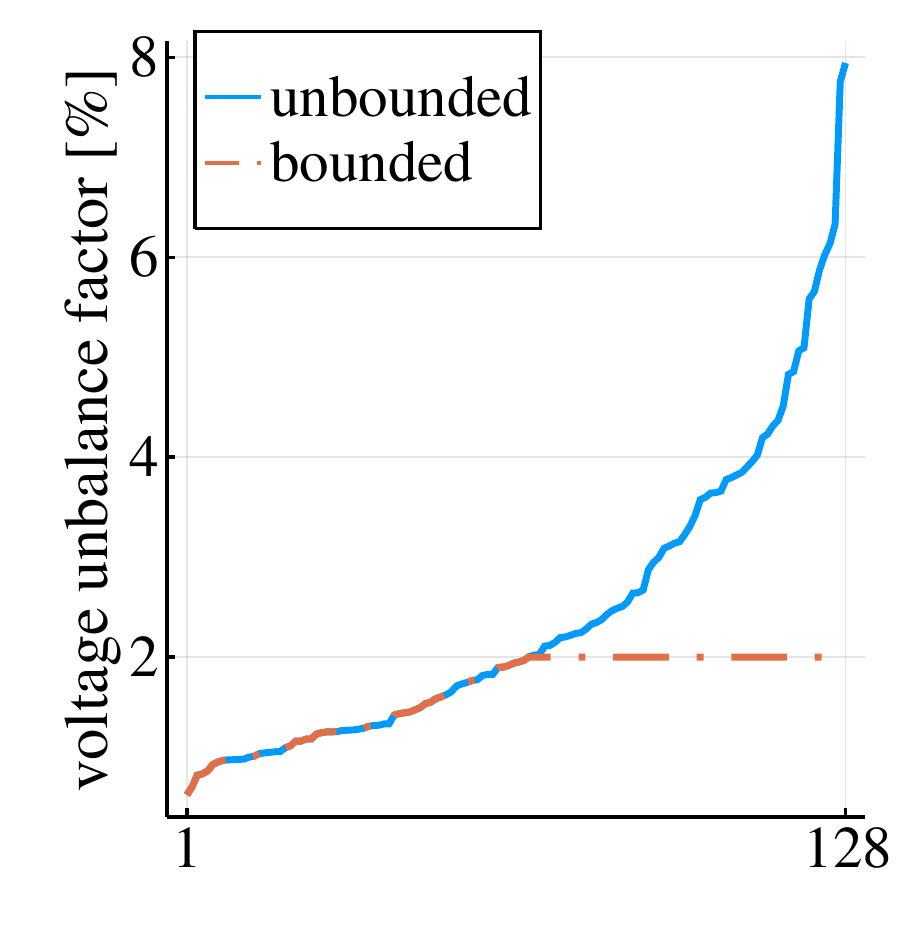}\\
	\hspace{2em}\footnotesize{network instance}\end{tabular}\label{fig:vufmax_10}}
	\hspace{0em}
	\subfloat[Kron-reduced]{\begin{tabular}[b]{c}\includegraphics[width=1.5in,trim=0 0.2in 0 0.1in,clip]{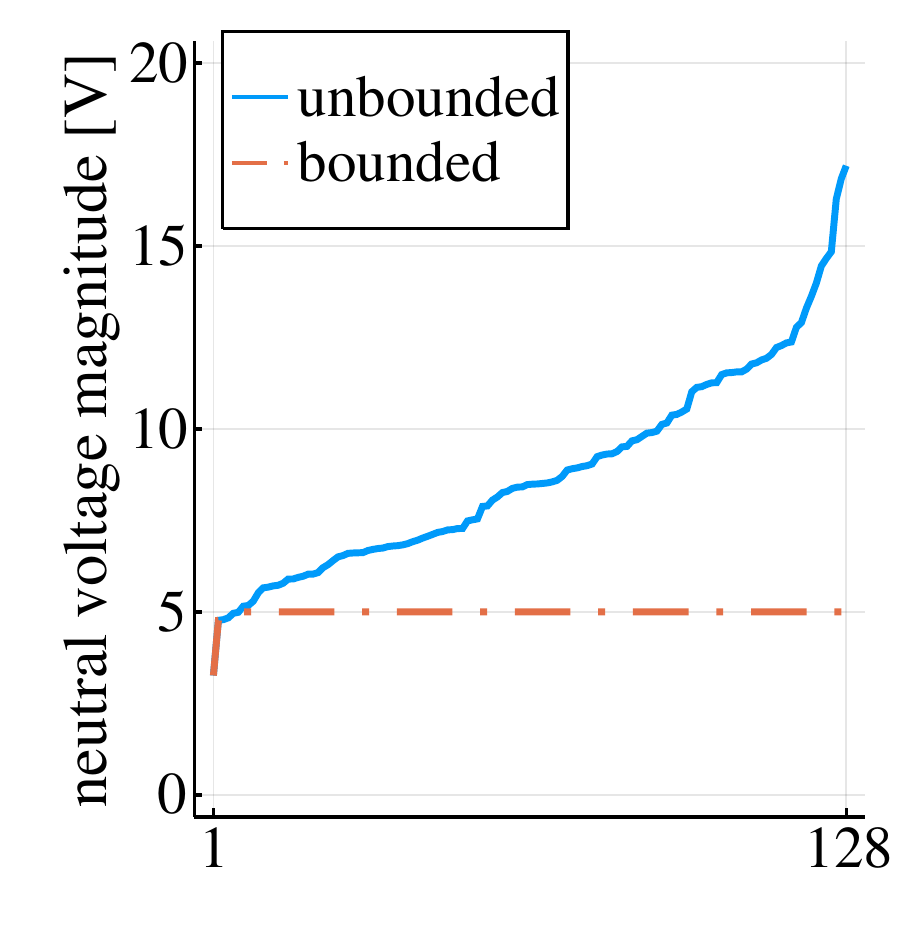}\\
	\hspace{2em}\footnotesize{network instance}\end{tabular}\label{fig:vmnmax_10}}
	\caption{When only bounding the phase-to-neutral voltages, the maximum VUF and neutral shift reach high values, which can be limited with explicit constraints (ordered by increasing value).}
	\label{fig:unbmax}
\end{figure}

Figure \ref{fig:unbmax} shows the highest VUF and neutral voltage magnitude for each instance when using the four-wire model and only imposing upper and lower bounds on the phase-to-neutral voltage magnitudes.
Figure \ref{fig:unbmax} also shows the resulting VUF and neutral shift when these are explicitly constrained to 2\,\% and 5\,V respectively.

\subsubsection{Solve time comparison}
\begin{figure}[tbh]
	\centering
	\includegraphics[width=3.5in]{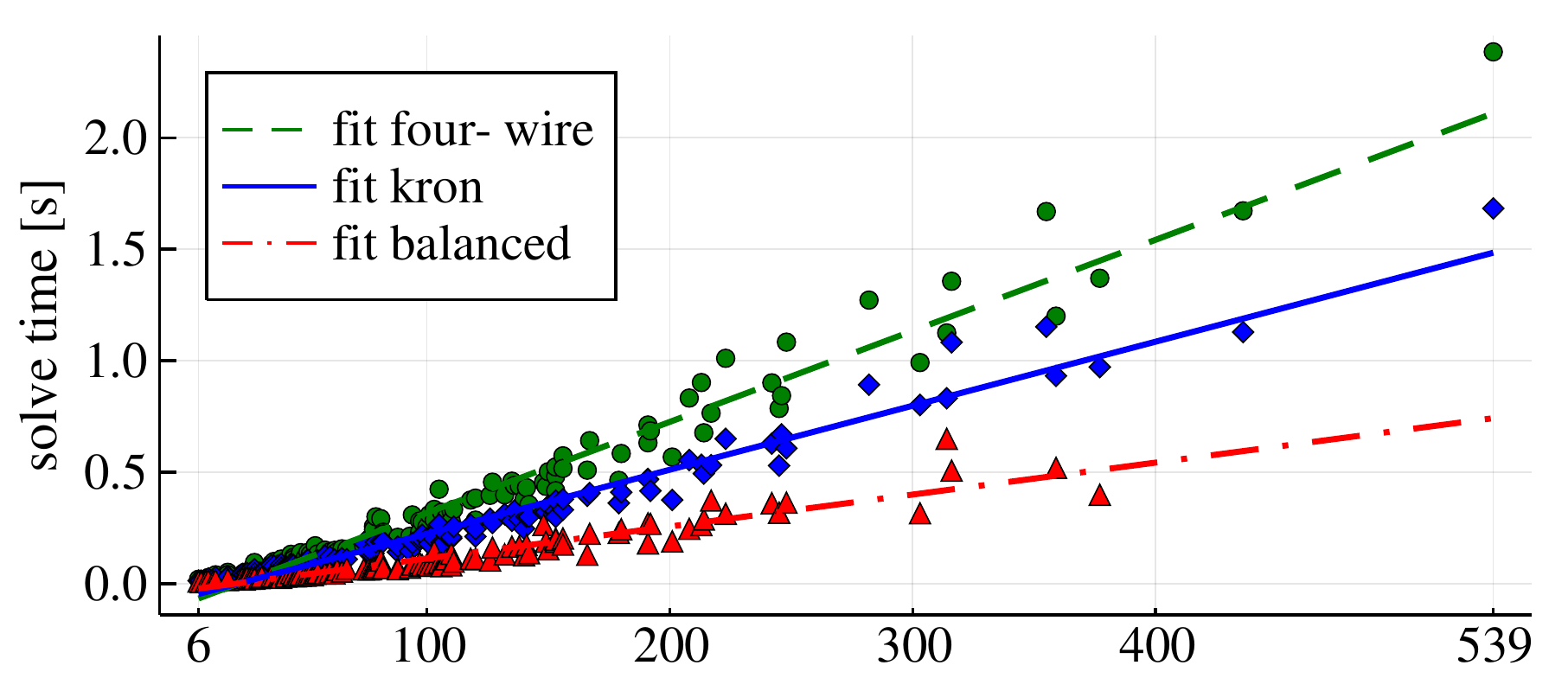}\\[-0.6em]
	\footnotesize{network size (number of buses)}
	\caption{Each marker shows the solve time for one of the 128 instances, together with the network size of that instance. This was done for the four-wire (circles), Kron-reduced (diamonds) and balanced (triangles) model. The lines show a least-square linear regression for each model respectively. For the balanced model, outliers at (246, 2.89s), (282, 7.06s) and (355, 6.63s) were left out.}
	\label{fig:solve_time_models}
\end{figure}

Figure \ref{fig:solve_time_models} shows how the solve time increases with the network size, for each model type. The solve time increases approximately linearly for all model types, up to networks with about 500 buses. The four-wire model takes about 40\,\% longer to solve than a Kron-reduced model. This corresponds approximately to a $\times$4/3 increase in the number of variables, as the four-wire model typically uses 4 variables for every 3 that the Kron-reduced model introduces.

%

\subsection{Power-voltage formulation}

\subsubsection{Solve time comparison}
\begin{figure}[tbh]
	\centering
	\includegraphics[width=3.5in]{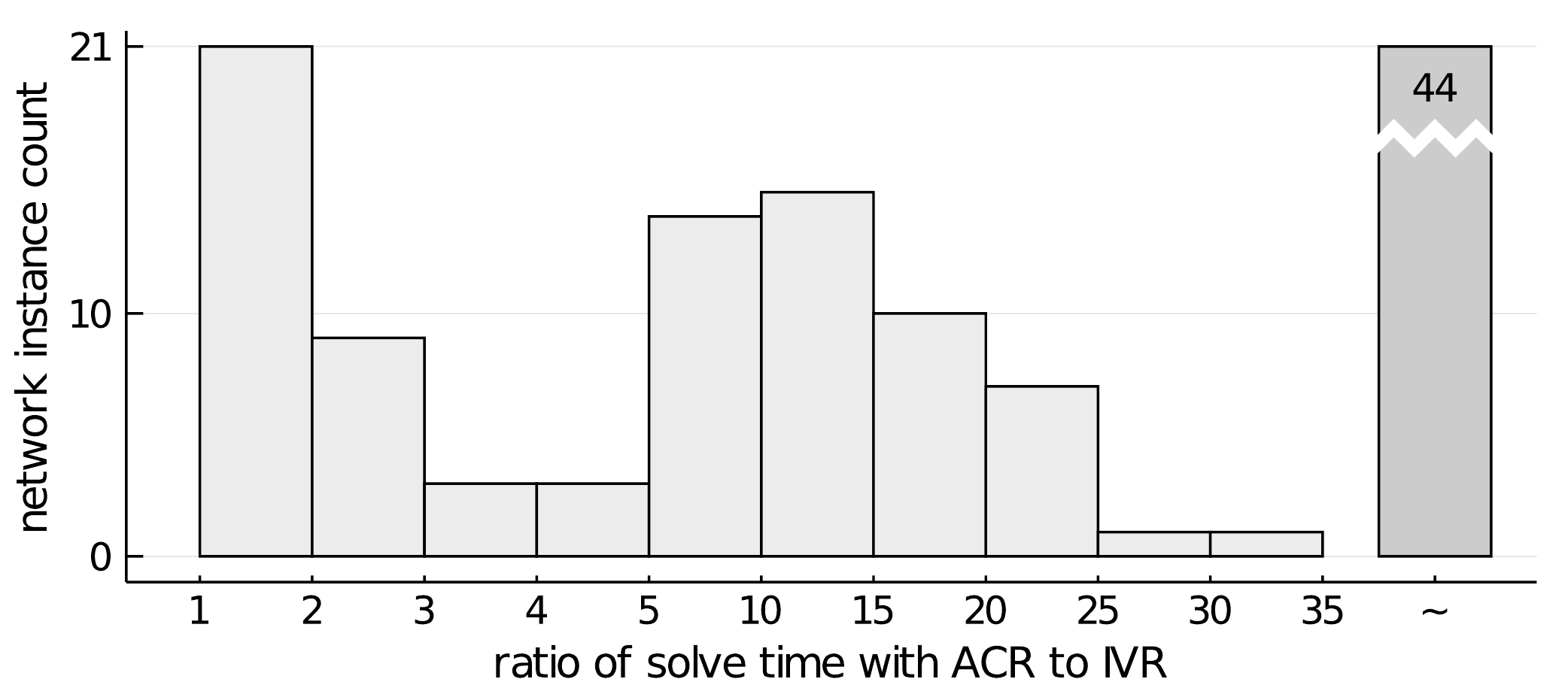}
	\caption{The power-voltage formulation takes a lot longer to solve for the four-wire model across all instances; 44 instances did not even converge within 500 iterations.}
	\label{fig:slowdown_acrivr}
\end{figure}

Figure \ref{fig:slowdown_acrivr} compares the solve time with the ACR  to IVR, for each of the 128 network instances. Both were initialized with the same, flat, balanced voltage profile. Unlike IVR, ACR did not even converge within 500 iterations for 35\,\% of the instances, and takes several times longer to solve.

\subsubsection{Constraint violation}
The authors argue that this is a consequence of the equivalence between KCL and the power balance, which breaks down when the nodal voltage goes to zero,
\begin{align}
	&\begin{aligned} &\hspace{1.5em}\text{KCL}\\&\sum \cscal{I} = 0\end{aligned}
	&&\begin{aligned} &\;\;\Rightarrow\\&\xLeftarrow[\cscal{U}\neq0]{}\end{aligned}
	&&\begin{aligned} &\hspace{-0.8em}\text{power balance}\\&\sum\cscal{S} = 0\end{aligned}.
\end{align}

\begin{figure}[tbh]
	\centering
	\newcommand{\vgfigheight}{0.4in}
	\subfloat[]{\includegraphics[height=\vgfigheight]{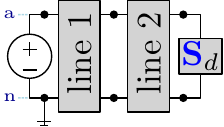}%
	\label{fig:vg_a}}
	\hspace{0.8em}
	\subfloat[]{\includegraphics[height=\vgfigheight]{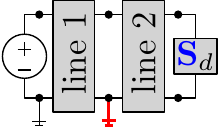}%
	\label{fig:vg_b}}
	\hspace{0.8em}
	\subfloat[]{\includegraphics[height=\vgfigheight]{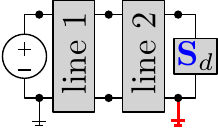}%
	\label{fig:vg_c}}
	\hspace{0.8em}
	\subfloat[]{\includegraphics[height=\vgfigheight]{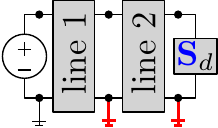}%
	\label{fig:vg_d}}
	\caption{Figure (a) shows the topology of a single-phase, 3-bus network as an example. The ACR formulation for (a) also allows all solutions to the networks (b-d), which are incorrect, non-physical solutions, introduced by the lifting of KCL to power variables.}
	\label{fig:vg}
\end{figure}

The power variables, across all components, are linked to the other variables by constraints of the form $\cscal{S}=\cscal{U}(\cscal{I})^*$. When $\cscal{U}=0$, $\cscal{S}=0$ , even when $\cscal{I}\neq 0$. Therefore, $\sum \cscal{S}=0$ is satisfied, whilst $\sum \cscal{I} = 0$ is no longer enforced. Clearly, all IVR-feasible solutions are feasible in ACR, but not the other way around. Therefore, the ACR formulation is a relaxation of IVR, which allows a non-physical grounding of any of the nodes. Figure \ref{fig:vg} summarizes this graphically.

These additional, incorrect solutions can be excluded by putting a lower bound on the voltage magnitude. 
In a four-wire network model, when a neutral terminal is not grounded, the voltage  will have a relatively small, non-zero magnitude. 
Since it can be arbitrarily close to zero depending on the network loading, it is very difficult to impose valid lower bounds a priori. 
The phase voltage variables however, will rarely drop below e.g. 0.5 pu; often, such bounds are even included in the problem specification. 
Because of this, this problem does not appear when using ACR to solve Kron-reduced and balanced network models.

\begin{figure}[tbh]
	\centering
	\includegraphics[width=3.5in]{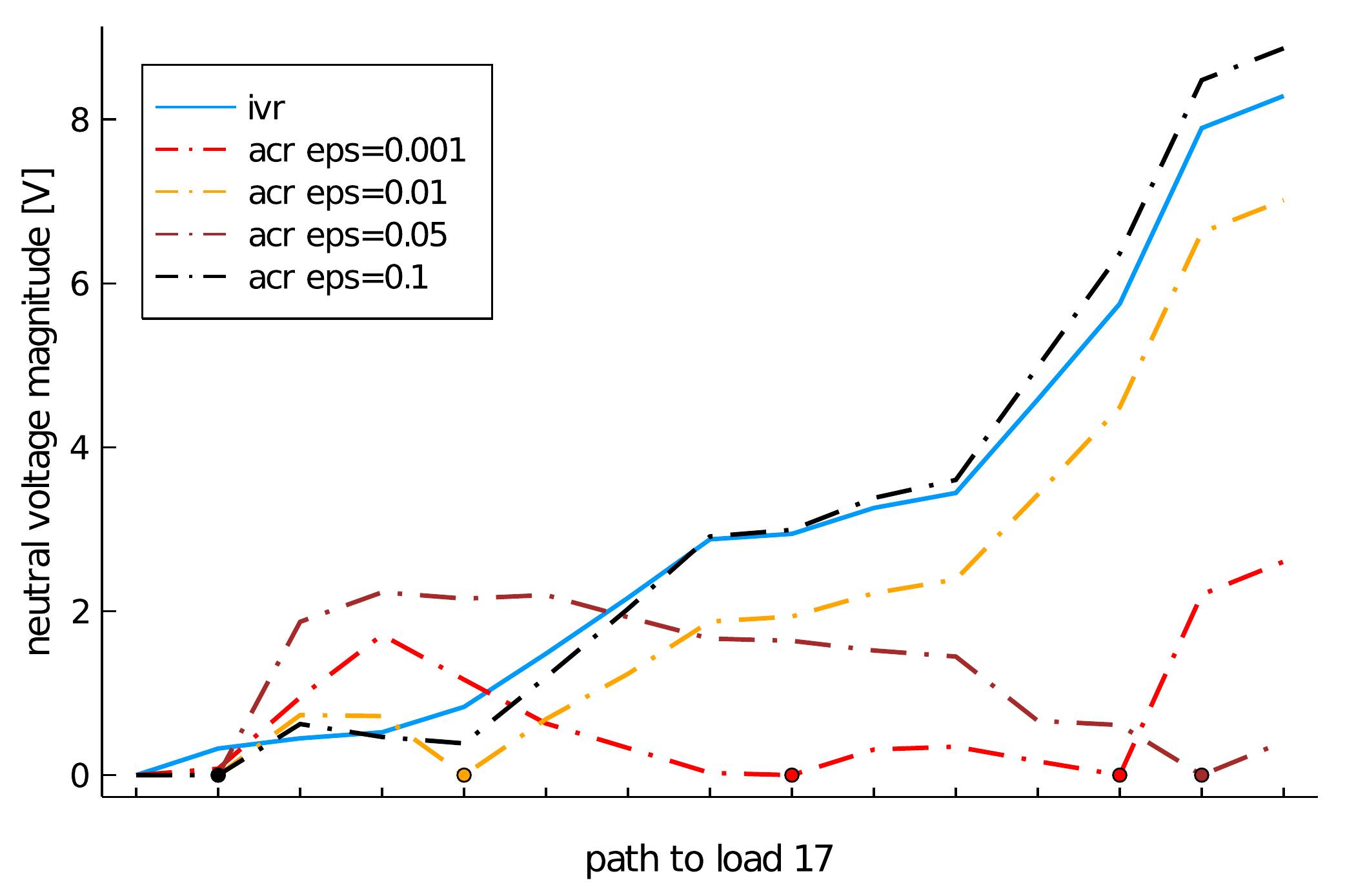}
	\caption{This figure shows the neutral voltage magnitude along the path from the substation to load 17, for the ENWL network instance (2,5), for solutions obtained with the IVR and ACR formulation. The dots indicate buses where the voltage magnitude is nearly zero ($<$1E-8).}
	\label{fig:vmn_acrivr}
\end{figure}

The incorrect solutions introduced by ACR, are equivalent to grounding the neutral at one or many buses.
Since doing this typically improves the network transfer capacity, it can be expected that this is also happens when maximizing the DG output in the case study.
Figure \ref{fig:vmn_acrivr} shows that this is exactly what happens.
The neutral voltage variables are initialized at $\epsilon+j0$ pu. 
The neutral voltage profile obtained with IVR, is only zero at the secondary of the transformer, where the neutral is actually grounded. 
Moving towards the load throughout the network, the neutral voltage steadily increases. Independent of $\epsilon$, the objective value always converged to -2.743.

ACR has many solutions, most of which correspond to solutions which are IVR infeasible. 
The result is very sensitive to the value of $\epsilon$, leading to different profiles for the neutral voltages with non-physical groundings at intermediate buses.
The objective value ranges from -2.780 to -2.757. 
Even with an initialization as high as 0.1\,pu, the problem persists. 
The authors argue that this is an inherent disadvantage of lifting KCL to power variables for four-wire networks, potentially leading to wrong solutions and unstable solver runs.

%
%

\begin{figure}[tbh]
	\centering
	\includegraphics[width=2.5in]{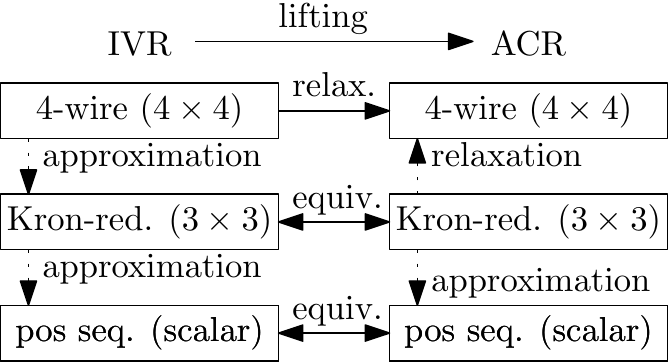}
	\caption{Summarizing the properties of the lifted and approximated formulations.}
	\label{fig:4wrelax}
\end{figure}

\subsubsection{Comparison of optimality}
Fig. \ref{fig:4wrelax} summarizes how the various model approximations and formulations are related. 
Reformulating KCL in terms of power flow variables introduces additional, non-physical solutions if the voltage magnitude is not constrained below. 
This is a relaxation, and allows non-physical solutions which can have a lower objective value. 
Also, note that these non-physical solutions can be interpreted as grounding any combination of neutral terminals. 
Therefore, the ACR four-wire formulation also includes all solutions of the ACR Kron-reduced formulation, i.e. when all neutral terminals are grounded. 

\begin{table}[]
    \centering
    \caption{Objective values for solutions obtained for a selection of ENWL test cases.  Note that `=' denotes that the ACR formulation yielded the same objective value as the IVR one.}
    \label{tab:obj_vals}
    \begin{tabular}{cllllll}\toprule
            & \multicolumn{2}{c}{four-wire} & \multicolumn{2}{c}{Kron-reduced} & \multicolumn{2}{c}{balanced}\\\cmidrule(lr){2-3}\cmidrule(lr){4-5}\cmidrule(lr){6-7}
        n-f & IVR & ACR & IVR & ACR & IVR & ACR\\
        \midrule
        04-2 & -1.5384 & -1.6700 & -1.6589 & = & -2.6392 & =\\
        23-3 & -0.7868 & -1.1020 & -1.1038 & = & -5.3645 & =\\
        10-5 & -2.1162 & -2.6966 & -2.7441 & = & -5.3718 & =\\
        11-4 & -6.3137 & -6.6470 & -6.7165 & = & -8.5624 & =\\
        09-3 & -10.022 & -11.452 & -11.984 & = & -17.529 & =\\
        18-5 & -4.1726 & -4.4614 & -4.6412 & = & -7.3296 & =\\
        05-2 & -2.4577 & -2.5544 & -2.6838 & = & -3.6570 & =\\
        11-1 & -5.8191 & -6.0016 & -6.4291 & = & -8.6007 & =\\
        23-1 & -4.6631 & -4.7166 & -4.8993 & = & -6.0126 & =\\
        12-2 & -3.9818 & -3.9793 & -4.2093 & = & -4.9350 & =\\
        19-4 & -3.3815 & -3.4719 & -3.3767 & = & -5.4268 & =\\\bottomrule
    \end{tabular}
\end{table}

Table \ref{tab:obj_vals} shows the objective value for a selection of ENWL test cases, obtained for all formulations shown in Fig.~a\ref{fig:4wrelax}. 
The objective value of the ACR four-wire solution is usually in-between the IVR four-wire and ACR Kron-reduced ones, but exceptions in both directions occur. 
Since ACR is a relaxation of IVR for the four-wire model, its global optimum should have a lower objective value. 
However, these optimization problems are non-linear and are only solved to \emph{local} optimality by \textsc{Ipopt}. 
Fig.~\ref{fig:vmn_acrivr} already illustrated how sensitive the obtained solution is to changes in the variable initialization. 
It is evident that the solution is not globally optimal for the entries where the objective value for the ACR four-wire solution is higher than the ACR Kron-reduced one, because it is a relaxation.

\section{Conclusions}
This article developed two quadratic formulations for unbalanced OPF, one in current-voltage variables and another one in power-voltage. Both forms include transformers, lines, shunts, loads and generators. The optimization models are agnostic of the number of phases, and can therefore model single-phase equivalent, Kron-reduced and four-wire networks.

The proposed current-voltage formulation is especially well-suited for four-wire networks.
It is shown that using  rectangular power-voltage variables  leads to problems, previously unobserved for positive-sequence and Kron-reduced representations.
When the neutral voltage is (close to) 0, the equivalence between KCL in current and power  breaks down.

A study across 128 low-voltage networks illustrated the effect of the network model detail on the set points obtained for the distributed generators. The Kron-reduced model tends to overestimate the network capacity, and caused the phase-to-neutral voltage magnitudes to exceed their upper bound (limited 10\,\% rise) by about 1-4\,\%. The solve time scales linearly with the number of variables, and the four-wire model takes therefore about 30\,\% longer to solve.


\section*{Acknowledgement}
S. Claeys carried out this research whilst holding a doctoral (PhD) strategic basic research grant (1S82518N) of Research Foundation - Flanders (FWO).
Andrew Urquhart and Murray Thomson at Loughborough University kindly provided the detailed line impedance data used in the case study.
We also thank David Fobes and Carleton Coffrin for their assistance  in integrating the proposed formulation into the open-source software package  \textsc{PowerModelsDistribution.jl} \cite{FOBES2020106664}.

\bibliographystyle{IEEEtran}

\end{document}